\documentclass[amsthm,10pt]{elsart}
\usepackage{amssymb}
\usepackage{amsmath,amsfonts}

\long\def\SKIPTEX#1\ENDSKIPTEX{}

\newtheoremstyle{usualtheorem}
  {\parskip}
  {}
  {\itshape}
  {}
  {\bfseries}
  {.}
  {.5em}
  {\thmnumber{#2}\thmname{ #1}\thmnote{ (#3)}}
\theoremstyle{usualtheorem}

\newcommand{\iffdef}{\rightleftharpoons}
\newcommand{\liff}{\longleftrightarrow}
\newcommand{\impl}{\longrightarrow}
\renewcommand{\lor}{\,\vee\,}
\renewcommand{\land}{\,\mathchar"2026\,}

\newcommand{\st}{\operatorname{st}}

\newcommand{\qsuper}[1]{^{\raisebox{.65\smallskipamount}{\scriptsize $#1$}}}

\def\rst_#1{\operatorname{st_#1}}
\newcommand{\srst}[1]{\qsuper{\st_#1}}

\newcommand{\sint}{\qsuper{\intr}}
\newcommand{\spart}[1]{{\vphantom{#1}}^\circ \kern-.1667em #1}
\newcommand{\stdn}[1]{{\vphantom{#1}}^{\mathrm s} \kern-.1667em #1}
\def\mon_#1(#2){\mu_{#1}(\,#2\,)}
\def\mons(#1){\mu(\,#1\,)}

\newcommand{\bq}{\;}
\newcommand{\All}{\boldsymbol\forall}
\newcommand{\all}{\bq\All}
\newcommand{\Ex}{\boldsymbol\exists}
\newcommand{\ex}{\bq\Ex}

\def\Allrst_#1{\All\srst{#1}}
\def\allrst_#1{\bq\Allrst_#1}
\def\Exrst_#1{\Ex\srst{#1}}
\def\exrst_#1{\bq\Exrst_#1}

\newcommand{\f}{\thickspace\thinspace}
\newcommand{\fb}{\medspace}

\newcommand{\Fl}[1]{\f#1\fb}
\newcommand{\fr}[1]{\fb#1}
\newcommand{\Fr}{\fr}

\newcommand{\Fbigl}[1]{\f\bigl#1\fb}
\newcommand{\fbigr}{\fb\bigr}
\newcommand{\Fbigr}{\fbigr}

\newcommand{\coll}[1]{\{\medspace #1 \medspace\}}
\newcommand{\scoll}[1]{{\vphantom{\{}}^{\rm s} \kern-.1667em \{\medspace #1 \medspace\}}
\renewcommand{\ss}{\subseteq}

\newcommand{\emp}{\emptyset}
\newcommand{\<}{\langle}
\renewcommand{\>}{\rangle}
\newcommand{\dom}{\operatorname{dom}}

\newcommand{\Fnc}{\operatorname{Fnc}}

\newcommand{\NCT}{{\bf NCT}}

\newcommand{\BST}{{\bf BST}}
\newcommand{\ZFC}{{\bf ZFC}}
\newcommand{\IST}{{\bf IST}}

\newcommand{\stan}{{\ifmmode \mathbb S\else $\mathbb S$\ \fi}}
\newcommand{\E}{{\ifmmode \mathbb E\else $\mathbb E$\ \fi}}

\newcommand{\Sms}{\operatorname{Sms}}
\def\Sat_#1{\operatorname{Sat_#1}}

\newcommand{\Set}{{\mathrm{Set}}}

\newtheorem{theorem}{Theorem}[section]
\newtheorem{proposition}[theorem]{Proposition}
\newtheorem{corollary}[theorem]{Corollary}

\newcounter{proofitem}

\def\vc#1*#2{#1_1,\ldots,#1_{#2}}
\def\vci#1#2*#3{#1_{#2},\ldots,#1_{#3}}
\newcommand{\prettytwo}[2]{%
$#1$\\
\hbox to .8\textwidth{\hfil$#2$}%
}

\newtheorem{remark}[theorem]{Remark}

\newtheoremstyle{fact}
  {\parskip}
  {0pt}
  {}
  {}
  {\bfseries}
  {}
  {.5em}
  {\thmnumber{#2}\thmnote{ #3}}
\theoremstyle{fact}

\newtheorem{fakt}[theorem]{}
\newcommand{\FACT}{\begin{fakt}}
\newcommand{\EFACT}{\end{fakt}}

\renewcommand{\emph}[1]{{\sl #1}} 

\newcommand{\THS}{{\ifmmode{\mathsf{THS}}\else{\sf THS}\fi}}
\newcommand{\TFS}{{\ifmmode{\mathsf{TFS}}\else{\sf TFS}\fi}}
\newcommand{\THSZ}{{\ifmmode \mathsf{THS}_0 \else $\mathsf{THS}_0$ \fi}}
\newcommand{\ZF}{{\ifmmode \mathsf{ZF}%
        \else $\mathsf{ZF}$ \fi}}
\newcommand{\ZFfin}{{\ifmmode \mathsf{ZF}^{\mathrm{fin}}%
        \else $\mathsf{ZF}^{\mathrm{fin}}$ \fi}}
\newcommand{\ZC}{{\ifmmode{\mathsf{ZC}}\else{\sf ZC}\fi}}
\newcommand{\ZCm}{{\ifmmode \mathsf{ZC}^{-}%
        \else $\mathsf{ZC}^{-}$ \fi}}

\renewcommand{\NCT}{{\sf NCT}}
\newcommand{\NGB}{{\sf NBG}}
\newcommand{\AST}{{\sf AST}}
\newcommand{\Vopenka}{Vop\v enka}
\newcommand{\opname}[1]{\operatorname{\mathsf{#1}}}
\renewcommand{\Set}{\opname{Set}}
\renewcommand{\st}{\opname{st}}
\renewcommand{\Sms}{\opname{Sms}}
\renewcommand{\Fnc}{\opname{Fnc}}
\renewcommand{\dom}{\opname{dom}}
\renewcommand{\phi}{\varphi}
\renewcommand{\emp}{\varnothing}
\newcommand{\smal}{\opname{small}}
\newcommand{\thin}{\opname{thin}}
\renewcommand{\inf}{\opname{inf}}
\newcommand{\bH}{\mathbb{H}}

\newcommand{\nat}{\mathbb{N}}
\newcommand{\snat}{\mathbb{SN}}
\newcommand{\Z}{\mathbb{Z}}
\renewcommand{\sint}{\mathbb{SZ}}
\newcommand{\R}{\mathbb{R}}
\newcommand{\sreal}{\spart{\mathbb{R}}}
\newcommand{\N}{\mathbb{N}}

\newcommand{\Q}{\mathbb{Q}}
\renewcommand{\a}{\alpha}
\renewcommand{\b}{\beta}
\renewcommand{\d}{\delta}

\newcommand{\s}{\sigma}
\newcommand{\srat}{\mathbb{SQ}}
\newcommand{\K}{\opname{ac}}

\newcommand{\cP}{\mathcal{P}}
\newcommand{\cX}{\mathcal{X}}
\newcommand{\cY}{\mathcal{Y}}

\newcommand{\restr}{\upharpoonright}

\newcommand{\Allthin}{\All\qsuper{\thin}}
\newcommand{\Exthin}{\Ex\qsuper{\thin}}

\newcommand{\exthin}{\ex\qsuper{\thin}}

\newcommand{\Allinf}{\All\qsuper{\inf}}

\newcommand{\ax}[1]{{\bf #1}}

\newcommand{\oneaxiom}[2]{\[\begin{tabular*}{\columnwidth}{|||lp{0.8\columnwidth}}{\bf
#1} & \hfil $#2$ \hfil\\\hline\end{tabular*}\]}

\newcommand{\defterm}[1]{\underline{\emph{#1}}}
\newcommand*{\eqdef}{\stackrel{\textup{def}}{=}}
\newcommand{\SL}{\opname{SL}}
\newcommand{\itemref}[1]{(\ref{#1})}
\newcommand{\Liff}{\Longleftrightarrow}
\begin{document}
\begin{frontmatter}
\author{P.~V.~Andreev \thanksref{RFBR}}
\address{1106-563 Zelenograd, Moscow, 124460, Russia}
\ead{petr@sherwood.ru}
\author{E.~I.~Gordon}
\address{
Department of Mathematics and Computer Science,
Eastern Illinois University,
600 Lincoln Avenue,
Charleston, IL 61920-3099,
USA
}
\ead{cfyig@eiu.edu}
\title{A Theory of Hyperfinite Sets}
\thanks[RFBR]{Partially supported by RFBR, Russia, grant N~03-01-00757.}
\begin{abstract}
We develop an axiomatic set theory --- the Theory of Hyperfinite Sets \THS,
which is based on the idea of existence of proper subclasses of big finite
sets. We demonstrate how theorems of classical continuous mathematics can be
transfered to {\THS}, prove consistency of \THS{} and present some
applications. 
\end{abstract}
\end{frontmatter}
\section*{Introduction}

Many applications of nonstandard analysis are based on
the simulation of infinite structures by hyperfinite ones. When
translated into the language of standard mathematics such
simulation means an approximation of infinite structures by finite
ones. Thus, nonstandard analysis provides us with a machinery
that allows to obtain new results about infinite structures using
such approximations and corresponding results about finite
structures. The latter are often much easier to obtain. This
approach is implemented in the famous monograph \cite{Nel} for the
construction of probability theory on infinite probability
spaces. In the monograph \cite{Gor} it was shown how this approach
can be used for systematic construction of harmonic analysis
on locally compact abelian groups starting from harmonic
analysis on finite abelian groups.

The results obtained on this way allow to look at this approach from another point of view.

According to this point of view Mathematics should be developed on
the base of the hypothesis that all sets are finite (some kind of
the ancient Greeks' atomism).

The historically first approach due to A. Yessenin-Volpin \cite{YV1}, \cite{YV2}
to develop this idea on the base of modern logic is called ultraintuitionism. It assumes the
existence of the maximal natural number.  This approach stimulated
investigations of the notion of feasible numbers. The first
mathematically rigorous formalization of the notion of feasibility
of natural numbers was introduced by R. Parikh \cite{Par}. Many
papers develop R. Parikh's approach as well as some other
approaches to the notion of feasibility (see e.g. \cite {Gu},
\cite{Dr}, \cite{Sa}). We do not discuss them here. A very
interesting discussion of correlation between the Real Analysis
and the Discrete Analyis is contained in \cite{Ze}. The main idea
of this paper is as follows: "Continuous analysis and geometry are
just degenerate approximations to the discrete world... While
discrete analysis is conceptually simpler ... than continuous
analysis, technically it is ususally much more difficult. Granted,
real geometry and analysis were necessary simplifications to
enable humans to make progress in science and mathematics....". In
some sense, our paper together with the paper \cite{gg}
contributes to this idea.

In this paper we develop an axiomatic theory of finite sets, which
we call the Theory of Hyperfinite Sets (\THS), by the reasons
explained below. Similarly to Kelley-Morse theory or von Neumann - Bernays -
G\"odel  theory
(\NGB), {\THS}  is a theory of classes in the $\in$--language,
where sets are defined as elements of classes. The universe of
sets satisfies all the axioms of \ZF$^{\textup{fin}}$ -- the
theory obtained by replacing in {\bf ZF} the axiom of infinity by
its negation and adding a suitable form of regularity (for
instance, an axiom saying that every set has a transitive closure;
see \cite{Soch}).

However, the properties of classes differ essentially from those
of \NGB. For example, the Separation Axiom fails in \THS: there
exist sets that contain proper subclasses (subclasses that are not
subsets). The reason why we need to include the last statement in
our theory is that we want to consider in {\THS} such properties
as feasibility discussed above. Indeed, let $F(x)$ be the
statement "$x$ is a feasible number" and $N$ be a non-feasible
number. Then the set $A=\{x\leq N\ |\ F(x)\}$ satisfies the
following inconsistent conditions: 1) $0\in A$, 2) $\forall x(x\in
A\longrightarrow x+1\in A)$, 3) $N\notin A$. The only way to avoid
this paradox, if one wants to keep the induction principle for
sets, is to assume that $A$ is not a set, and thus the separation
axiom fails for the finite set $\{0,1,\dots,N\}$.

The paradox discussed in the previous paragraph is a version of
the well-known paradox about a pile of sand, due to Eubilides, IV century
B.C.: since one grain of
sand is not a pile and if $n$ grains of sand do not form a pile of
sand, then $n+1$ grains do not form a pile of sand also, then how
can we get a pile of sand? The paradoxes of these type can not be
considered in the framework of classical set theory since the
objects, like a pile of sand, have a very vague description and,
thus, cannot be considered as any objects of classical
mathematics, i.e. as sets. On the other hand there are many 
examples that show that such notions arise very naturally in
mathematics (see, e.g., the example concerning the feasibility
above). The first mathematician who realized the importance of the
notions  of this type was P.\Vopenka. In \cite{Vop} he introduced
the first theory of finite sets, the Alternative Set Theory
(\AST), where the existence of finite sets containing subclasses
that are not sets was postulated. Such subclasses of sets are
called {\it semisets}.

The main defect of P.~\Vopenka's approach is the opposition of his
theory to classical mathematics. As it was mentioned above (see
the quotation from \cite{Ze}) the advantage of the continuous
mathematics in comparison with the discrete one is its simplicity 
that often allows to solve problems concerning discrete objects.

\THS{} introduced here is also based on the idea of existence
of proper subclasses of big finite sets. Finite sets that contain
proper subclasses are called hyperfinite sets. This term is
borrowed from nonstandard analysis. The primary model for \THS{}
is the collection of all subclasses of the set of hereditarily
finite sets in the Nonstandard Class Theory {\NCT} \cite{AG}. The
central notion of a {\emph thin class} is defined by a formulation
equivalent in \NCT{} to the definition of a class of standard size:
a class is thin if any subset of it does not contain proper
subclasses. Sets that do not contain proper subclasses are called
{\it small}. The
class of all small natural numbers is a thin class. It
coincides with the set $\omega$ in \ZF. Under our approach the
class of all small numbers can be interpreted as the class of
feasible numbers.

We prove that all results of classical mathematics that can be
formalized in Zermelo set theory can be proved for thin structures in \THS. 
This is a substantial difference between \THS{} and \AST. It allows to
formalize within \THS{} those proofs of theorems about finite sets
that use continuous mathematics and, hence, it is not necessary to
invent any new proofs for such theorems.

In the discrete world continuous objects have their place as well:
they originate from hyperfinite sets or their $\s$-subclasses as
quotient "sets" by some {\it indiscernibility} relation. An
indiscernibility relation $\rho$ is an equivalence relation that
is a $\pi$-class and satisfy some special condition (see 
section~\ref{indiscernibility}). A class is called a $\s$-class ($\pi$-class) if it can be
represented by the union (the intersection) of a thin class. We
prove that there exists a thin class of representatives of all
$\rho$-equivalence classes, which represent the quotient "set" (more
exactly, the quotient system of classes) by $\rho$.

For example, to obtain the field of reals $\R$ in \THS{} one
should consider a computer arithmetic implemented in an
idealized computer with a hyperfinite memory for simulation of the
field of reals.  It may be the usual computer arithmetic, based on
the representation of reals in the form with floating point. Let
$\<R;\oplus,\odot\>$ be this system. It is well-known that because
of the rounding off the operations $\oplus$ and $\odot$ are neither
associative, nor distributive. Let $R_b=\{ x\in R\ |\
\exists^{\textup{small}} n (|x|<n)\}$, where $\exists^{\textup{small}}
n$ means "there exists a small natural number $n$". Since the
class of all small natural numbers is a thin class, it is easy
to see that $R_b$ is a $\s$-class. We can interpret the elements
of $R_b$ as the computer numbers that are far enough from the
boundary of the computer's memory, so that doing computations with
these numbers one can never get overfilling of memory. Indeed,
it can be proved that the class $R_b$ is closed under the
operations $\oplus$ and $\odot$. The indiscernibility
 relation
$\rho$ is defined by the condition $x\rho y\Longleftrightarrow
\forall^{\textup{small}} n (|x-y|<\frac 1n)$. Obviously $\rho$ is a
$\pi$-class. It also has the natural interpretation: we identify
those numbers that differ on a number close enough to the computer
zero. It can be proved (in \THS) that the quotient system
$R_b/\rho$ is isomorphic to the field $\R$.

Certainly, there are many other systems, from which one can obtain
$\R$ in the way similar to one described in the previous
paragraph. The system based on representation of reals in the form
with floating point is discussed in details in \cite{ggh}. In this
paper we introduce a hyperfinite system that is a little bit
simpler and has some better properties - it is an abelian group
for addition. However, we proved that it is impossible to obtain
the field $\R$ from a hyperfinite system that is an associative
ring \cite{gg}. Similar facts hold also for many locally compact
non-commutative groups. It is shown in \cite{gg} that we cannot
find the hyperfinite groups that have approximate properties of
many important Lie groups such as $SO(3)$. The hyperfinite objects
with the best properties, from which all unimodular locally
compact groups can be constructed, are quasigroups (latin squares)
\cite{gg}. These facts demonstrate that the continuous world has
better properties than the discrete one. The theory \THS{}
introduced here allows to formalize not only all classical
mathematics, but also the statements about the connection between
the discrete world and its continuous approximation.

In \cite{BF} the Non Standard Regular Finite Set Theory was formulated by
S.~Baratella and R.~Ferro. Based on a countably saturated universe of 
hereditarily hyperfinite sets, NRFST contains a rich structure of external sets
\emph{over} it. We believe that a theory of classes of higher levels over
the universe of hyperfinite sets of {\THS} can be formulated in the pure
$\in$--language and simulated \emph{within} \THS.

We are very grateful to Karel Hrbacek, Vladimir Kanovei, 
Edward Nelson and Anton\'{\i}n Sochor for
interesting and helpful discussions on this work.

\section{Axioms}

\begin{remark}
In~\cite[\S 6]{AG} we announced a theory of hyperfinite sets \THS{}. 
The theory presented here is a result of
further development of the idea; the implementation is very different though 
and, we believe, much more interesting than the one described in~\cite{AG}.
\end{remark}

\begin{fakt}
\THS{} is a first-order theory. 
Semantical objects of the theory are \defterm{classes}.
Its language contains only one non-logical symbol --- the binary 
predicate symbol $\in$ of the membership.
\end{fakt}

\begin{fakt}
\defterm{Sets} are \emph{defined} as members of classes:
\[
    \Set(X) \iffdef \Ex Y \Fl( X \in Y \Fr).
\]

We accept the convention to use small letters for sets
and capital letters for classes.
\end{fakt}

\FACT
Formulas where all quantifiers range over set variables are called
\defterm{normal formulas}.
\EFACT
\FACT
\defterm{Set formulas} are normal formulas where
no class constants or variables occur.
\EFACT

{\sf Axiom of Extensionality}:
\oneaxiom{Ext}{\All X \all Y \Fbigl( \All x 
\Fl( x\in X \liff x\in Y \Fr) \impl X=Y \Fbigr).}

{\sf Axioms of class formation} (arbitrary formulae are allowed):
\oneaxiom{Class}{\All X_1,\ldots,X_n \ex Y \all y 
\Fl( y\in Y \liff \Phi(y,X_1,\ldots, X_n) \Fr).}

{\sf Axiom of set formation}:
\oneaxiom{Set}{ \Set \emp \land \All x \all y \f \Set(x\cup
\coll{y} ) .}

{\sf Axioms of induction and regularity} (only set formulae are allowed):
\oneaxiom{Ind}{\phi(\emp) \land \All x \All y \Fl( \phi(x) \land \phi(y) 
\impl \phi(x\cup\coll{y} ) \Fr) \impl \All x \f \phi(x).}

\begin{fakt}
The class of all sets is denoted by $\bH$.
\end{fakt}

\begin{fakt}
Subclasses of sets are called \defterm{semisets}.

Sets in \THS{} may contain proper subsemisets, 
i.e. subclasses which are not sets.
\end{fakt}

\FACT
A set is called \defterm{small} if it does not contain proper subsemisets:
\[
    \smal x \iffdef \All Y\ss x \Fl( \Set Y \Fr).
\]
\EFACT

\FACT
A class $X$ which is not small is called \defterm{infinitely large} or
simply \defterm{infinite} ($\inf X$).
\[
    \inf X \iffdef \Ex Y\ss X \Fl( \lnot \Set Y \Fr).
\]
\EFACT

\FACT
A class is called \defterm{thin} iff every sub\emph{set} of it is small:
\[
    \thin X \iffdef \All a\ss X \all C\ss a \Fl( \Set C \Fr).
\]

Thus, thin set is the same as small set.

\EFACT

In formulas we use quantifiers with superscripts ${\thin}$, $\smal$ and
$\inf$ in a natural way.

{\sf Axiom of Thin Semisets}:
\oneaxiom{Thin}{\All X \Fbigl( \thin X \impl \Ex x \Fl( X\ss x \Fr) \Fbigr).}

{\sf Axiom of Compactness}:
\oneaxiom{Comp}{\Allthin X \all u \Fbigl( u \ss \cup X
\impl \Ex x\ss X \Fl( u\ss\cup x \Fr) \Fbigr).}

{\sf Axiom of Exponentiation}:
\oneaxiom{Exp}{\Allthin X \exthin P \all y \ex p\in P 
\Fl( y\cap X = p\cap X \Fr).}

\begin{fakt}
We define the ordered pair $\<x,y\>\eqdef\coll{\coll{x},\coll{x,y}}$ and the operations
$\times$~(cartesian product), $\dom$ (domain), $"$ (image: $X"A=\coll{b:\Ex a\in A
\Fl( \<a,b\> \in X \Fr) }$) in the usual way. We also use $\Fnc F$ as a
shorthand for $\All x\all y\all z \Fl( \<x,y\> \in F \land \<x,z\>\in F
\impl y=z \Fr)$.
\end{fakt}

{\sf Axiom of Choice}:
\oneaxiom{Choice}{\All X \Fbigl( \thin \dom(X) \impl $\hfill\hfill\\
& \hfill $ \Ex F \Fl( \Fnc F \land \dom(F) = \dom(X) \land F\ss X \Fr) \Fbigr).}

\begin{fakt}
The class $\snat$ of \defterm{small natural
numbers} is defined as the smallest class which contains the empty set and is closed
under the von Neumann successor operation:
\[
    \snat = \coll{ x: \All N \Fbigl( \Fl[ \emp \in N \land  \All n\in N
        \Fl( n\cup \coll{n} \in N \Fr) \Fr] \impl x\in N \Fbigr) }.
\]
\end{fakt}

{\sf Axioms of Dependent Choices} (arbitrary formulae are allowed):
\oneaxiom{DC}{
\All X \ex Y \f \Phi(X,Y) \impl $\hfill\hfill\\
& $ \All X_0 \ex Z \Fbigl( Z"\coll{\emp}=X_0
\land \All n\in \snat \f \Phi( Z"\coll{n}, Z"\coll{ n \cup \coll{n} } )
\Fbigr).
}

\begin{fakt}
Similarly to small natural numbers, the class $\stan$ of all \defterm{standard sets} is defined as the smallest
class containing the empty set and closed
under the operation of adjoining one element:
\[
    \stan = \coll{ x: \All S \Fbigl( \Fl[ \emp \in S \land  \All a\in S \all b\in S
        \Fl( a\cup \coll{b} \in S \Fr) \Fr] \impl x\in S \Fbigr) }.
\]
\end{fakt}

{\sf Axioms of Transfer} (only set formulae are allowed):
\oneaxiom{T}{\All t_1\in\stan \cdots \all t_n\in\stan 
\Fl( \Ex x \f \phi(x,t_1,\ldots,t_n)$\hfill\hfill\\
& \hfill $ \impl \ex x\in\stan \f \phi(x,t_1,\ldots,t_n) \Fr).}

\begin{fakt}
We denote
\begin{eqnarray*}
\mathsf{TFS} & = &\ax{Ext} + \ax{Class} + \ax{Set} + \ax{Ind}\\
\THSZ & =  & \mathsf{TFS} + \ax{Thin} + \ax{Comp}\\
\THS & = & \THSZ + \ax{Exp} + \ax{Choice} + \ax{DC} + \ax{T}
\end{eqnarray*}
\end{fakt}

\begin{remark}
Axioms of \TFS{} are borrowed from \Vopenka's \AST.
\ax{Thin} and \ax{Comp} are true in \AST{} for countable classes.
See also~\ref{two-cardinalities-remark}, \ref{prolongation-remark} and
\ref{ast-zermelo}.
\end{remark}

\begin{remark}
The axioms of transfer are not as important in {\THS} 
as in other non-standard frameworks because standard sets are not the
primary object of investigation here. The main reason for including them
into the list of axioms is Theorem~\ref{finite-rich}.
\end{remark}


\section{Basic facts and notions}

\begin{fakt}
\defterm{Natural numbers} are defined the same way as ordinals are defined in \ZF:
they are transitive sets linearly ordered by the membership
relation. The class of all natural numbers is denoted by $\nat$. 
{\bf Ind} implies induction over $\nat$ for any set-formula $\phi$:
\[
    \Fl[ \phi(\emp) \land \All n\in\nat \Fl( \phi(n) \impl \phi(n\cup\coll{n}) \Fr) \Fr]
    \impl \All n\in\nat \f\phi(n).
\]
\end{fakt}

\FACT
\ax{Ind} implies also that for any set $x$ its \defterm{size}
$\sharp(x)$ is
uniquely defined as a natural number $k$ such that there is a set-bijection
from $x$ onto $k$.
\EFACT

\begin{theorem}\ \label{zffin}
\begin{enumerate}
\item
All axioms of $\TFS$ hold in $\stan$.
\item
The universe $\bH$ of all sets and the universe $\stan$ of
standard sets both satisfy the axioms of \ZFfin.
\end{enumerate}
\end{theorem}
\begin{proof}
It follows
from the definition of $\stan$ that the axioms of {\TFS} are true in $\stan$.
Sochor\cite{Soch} proved that \ZFfin is equivalent to the theory
with the axioms \ax{Set}, \ax{Ind} and extensionality for sets. 
\end{proof}

\begin{proposition}[\TFS]
There exists a bijective mapping $\K$ from the universe $\bH$ of all sets
onto the class $\nat$ of natural numbers, definable by a set formula and
such that $x\in y$ implies $\K(x)<\K(y)$ for all sets $x$ and $y$.
\end{proposition}
\begin{proof}
It can be proved in \ZFfin that the Ackermann encoding of finite sets
defined inductively by the conditions $\K(\emp) = 0$ and $\K(x)=\sum_{a\in x} 2^{\K(a)}$
is a total bijection.
\end{proof}

\begin{proposition}[\TFS]\ \begin{enumerate}
\item
    $\stan$ is a thin class;
\item
    $\stan$ coincides with the class of hereditarily small sets;
\item
    $\snat = \coll{n\in\nat:\smal n}=\stan\cap\nat$.
\end{enumerate}
\end{proposition}

Many properties of small sets and thin classes can be proved in {\TFS}
already.

\begin{proposition}[\sf{TFS}]\label{tfs:props}\ \begin{enumerate}
\item\label{thinz} $\smal x \land y\ss x \impl \smal y$; 
                    $\thin X \land Y\ss X \impl \thin Y$;
\item\label{thini} $\smal x \liff \smal\sharp(x)$;
\item\label{thinii} $\smal x \impl \smal F\restr x \land \smal F"x$, for any function $F$;
\item\label{thiniii} $\Allinf n\in \nat \ex x \Fl( X\ss x \land \sharp(x) \le n \Fr) 
    \impl \thin X$;
\item\label{thiniv} $\smal \cup x \impl \smal x$; $\inf x \impl \inf \cup x$;
\item\label{thinvi} $\thin X \impl \thin \coll{ y: y\ss X }$;
\item\label{thinvii} $[\f \thin X \land \All x\in X \f \thin Y"\{x\} \Fr]
    \impl \thin Y\restr X$;
\item\label{thiniiip} $\thin X\land \thin Y \impl 
    \thin X\cup Y \land \thin X\times Y$;
\end{enumerate}
\end{proposition}
\begin{proof}
\itemref{thini}. $\impl$ Assume $\inf\sharp(x)$. Let $f$ be some set-bijection from
$\sharp(x)$ onto $x$. Then the class $Y=\coll{f(n):n\in\snat}$ is a proper
semiset, since otherwise $\snat=f^{-1}"Y$ would be a set. Thus, $x$ is
infinite.

$\longleftarrow$ One should proceed by induction on $\sharp(x)$ over $\snat$. Due to
\ax{Set} adjoining one element to a small set gives a small set again.

\itemref{thinii}. We use \itemref{thini} and proceed by induction over $\sharp(x)$.

\itemref{thiniii}. Let $X$ be not thin. Then, by definition of thin class, there
exists an infinite subset $y\ss X$. Therefore, every superset $x\supseteq X$
cannot contain less than $\sharp(y)$ elements.

\itemref{thiniv}. Let $y=\cup x$ be a small set. Then
$\sharp(x) \le 2^{\sharp(y)}$. Since, by \itemref{thini}, $\sharp(y)$
is small, $\sharp(x)$ is small and $x$ is small.

\itemref{thinvi}. Denote $Y= \coll{ y: y\ss X }$ and assume $u\ss Y$ is
infinite. Then, according to \itemref{thiniv}, $\cup u$ is infinite as well. But
$\cup u \ss X$, in contradiction with the fact that $X$ is thin.

\itemref{thinvii}. Let the left hand side of the implication holds. Take any set
$a\ss Y\restr X$. Then $\dom a$ is small, since $\thin X\supseteq \dom a$.
By \itemref{thini} the numbers $p=\sharp(\dom a)$ and 
$q=\max \coll{\sharp(a"u):u\in\dom a}$ are small. Hence, 
$\sharp(a) \le p\cdot q$
is also small, and $a$ is small. This proves that $Y\restr X$ is thin. By
\itemref{thiniiip}, $Y"X$ is also thin.
\end{proof}

\begin{proposition}[\THSZ]\label{thin-equiv-defs}
The following statements are equivalent for any class $X$:
\begin{enumerate}
\item $X$ is a thin class (all subsets of $X$ are small);
\item $\Allinf n \in \nat \ex a \Fl( X\ss a \land \sharp(a) = n \Fr)$;
\item $\All Y\ss X \ex y \Fl( Y = y\cap X \Fr)$.
\end{enumerate}
\end{proposition}
\begin{proof}
Assume $X$ is thin, $Y\ss X$ and $n\in\nat$ is infinitely large. 
Taking into account item \itemref{thiniii} 
of Proposition~\ref{tfs:props}, it is enough to show that 
\begin{equation}\label{two-rabbits}
    \Ex y\ss x \Fl( y\cap X=Y \land \sharp(y) \le n).
\end{equation}
By \ax{Thin} $X\ss x$ for some infinite set $x$. Denote 
$s=\coll{y\ss x : \sharp(y) \le n}$, 
$D = \big\{ \coll{ y\in s: a\notin y }: a\in Y\big\}$,
$\bar{D} = \big\{ \coll{ y\in s: b\in y }: b\in X\setminus Y\big\}$. 
$D$ is thin since for every set $t\ss D$ there is a set 
$\cup\coll{ x\setminus \cup d : d\in t} \ss Y$ of the same size. Similarly, 
$\bar{D}$ is thin as well.
Suppose \itemref{two-rabbits} does not hold. Then 
$\cup (D\cup \bar{D}) \supseteq s$. Hence, by \ax{Comp},
$\cup t \supseteq s$ for some small $t\ss D\cup\bar{D}$ which is impossible because
$\sharp(t)<\sharp(x)$.
\end{proof}

As an immediate corollary we get the following proposition.
\begin{proposition}[\THSZ+\ax{Exp}]\label{powerclass-codable}
$\Allthin X \exthin P \all Y\ss X \ex p\in P 
\Fl( Y = p\cap X \Fr).$
\end{proposition}

\begin{proposition}[\THSZ]\label{thin:union}\ 
\begin{enumerate}
\item
    $\thin X\impl \thin(\dom X)$;
\item\label{thinviii} 
    $\thin X \impl \Fbigl[ \thin F"X \land \All y\ss F"X 
    \ex x\ss X \Fl( y = F"x \Fr) \Fbigr]$, for any function $F$;
\item
    $\thin X \impl \Sms \cup X$.
\end{enumerate}
\end{proposition}
\begin{proof}
(1) follows immediately from the previous proposition since 
$\sharp(\dom(x))<\sharp(x)$ for any $x$.

\itemref{thinviii}. Let $X$ be thin. It follows from
Proposition~\ref{tfs:props},\itemref{thinvii} that $F\restr X$ is
also thin. By 1 $F"X$ is thin. Using induction over $\snat$ on the cardinality of $y\ss F"X$ one
proves the existence of a set $x$ such that $F"x=y$.
\end{proof}

\begin{fakt}
A class is called \defterm{countable} iff it can be bijectively mapped onto
the class of small natural numbers.
\end{fakt}

\FACT \label{two-cardinalities-remark}
It follows from the previous proposition that the
theory \AST+\ax{Thin}+\ax{Comp}+"there exists an uncountable thin class" is inconsistent.
Indeed, in \AST{}, due to the axiom of two cardinalities saying that there
is a bijection between any two uncountable classes, an 
uncountable thin semiset can be bijectively mapped onto an
infinite set which is not thin.
\EFACT

\FACT
The property of being a thin infinite class behaves as a cardinality lying between the cardinalities
of small sets and those of infinite sets. This fact can be expressed in the
following way. 

We define \emph{inner cardinality} of a class:
\[
    \opname{ICard} X \eqdef \coll{n:\Ex x \Fl( x\ss X \land \sharp(x) = n+1 \Fr) }.
\]
Then for sets we have $\opname{ICard} x = \sharp(x)$, and 
infinite thin classes are exactly the classes $X$ such that 
$\opname{ICard} X = \snat$.
\EFACT

\begin{proposition}[TFS]\label{prolongation}
Axioms \ax{Thin} and \ax{Comp} together are equivalent to the following
statement:
\begin{description}
\item[\sf Prolongation principle:]
\item[\qquad]
    $\All\qsuper{\thin} X \Fbigl[ \All\qsuper{\smal} x\ss X \fb \phi(x) \impl \Ex y \Fl( \Set(y) \land
    X\ss y \land \phi(y) \Fr) \Fbigr]$
\item[\qquad]
    where $\phi$ is any set-formula with set-parameters.
\end{description}
\end{proposition}
\begin{remark}\label{prolongation-remark}The prolongation principle formulated here is a generalization
of prolongation axiom of {\AST}, which says that every countable function is a
subclass of a set-function.
\end{remark}

The next proposition lists counterparts of statements which became customary
tools in non-standard analysis. All of them are just special cases of the
prolongation principle formulated above.

\begin{proposition}[$\THSZ$]\label{saturation-extension}\ 
\begin{description}
\item[\sf Saturation:] 
    $\Allthin Y \Fl( \All y\ss Y \Fl( \cap y \neq \emp
    \Fr) \impl \cap Y \neq \emp \Fr)$;
\item[\sf Extension:] 
    $\All F \Fbigl( \thin F \land \Fnc F \impl \Ex f \Fl( \Fnc f
    \land F\ss f \Fr) \Fbigr)$;
\item[\sf Nelson's idealization principle:]
\item[\qquad]
    $\Allthin A \Fbigl( \all a_0\ss A \ex x \all a\in a_0 \f \phi(a,x) \impl
        \ex x \all a\in A \f \phi(a,x) \Fbigr)$,
\item[\qquad]
    for any set--formula $\phi$ with set parameters.
\end{description}
\end{proposition}

\FACT
As we said already in the introduction, the simplest and the most important
proper classes are $\sigma$--classes and $\pi$--classes which are defined in
\THS{} as follows: a \defterm{$\sigma$--class} is a union of a thin class
and a \defterm{$\pi$--class} is an intersection of a thin class.

Both $\pi$--classes and $\sigma$--classes are semisets
(see~Proposition~\ref{thin:union}).
\EFACT

\FACT \label{system-of-classes}
Together with sets and classes one can consider in \THS{} also
\defterm{systems of classes} defined as collections of classes satisfying a
certain formula and written as terms of the form
\[
    \coll{ X : \Phi(X) },
\]
where $\Phi$ is an arbitrary formula with some class-- or set--parameters.
\EFACT

\FACT
A system of classes $\mathcal{X} = \coll{ X : \Phi(X) }$ is called \defterm{codable} if
there exists a class $C$ such that
\[
\mathcal{X} = \coll{ C"\coll{d} : d\in \dom(C) }.
\]

Such coding by a class $C$ is called \defterm{extensional} iff $C"\coll{d}
\neq C"\coll{d'}$ for distinct elements $d,d'\in \dom(C)$.

If a system $\cX$ is coded by a class $C$ one can use quantification over
subsystems of $\cX$. For instance, $\All \cY\ss\cX\f \Phi(\cY)$ can be
interpreted as $\All D\ss\dom C \f \Phi( \coll{ C"\coll{d} : d\in D } )$.

\EFACT
\begin{fakt}
If there exists a coding $C$ such that $\dom(C)$ is a thin class, 
the system $\mathcal{X}$ is called a \defterm{thin system of classes}. 

If $\mathcal{X}$ is a thin system we can always assume without loss of 
generality, due to axiom of Choice, that a given coding of $\cX$ is
extensional. 

If, furthermore, $\cX$ is a system of thin classes we can speak of 
systems of systems of classes and so
on, using an appropriate encoding for higher levels. Such an encoding can
always be chosen to be a thin class, due to the axiom of exponentiation.

\end{fakt}

\begin{proposition}[\THS]\label{thin-wo}\ 
\begin{enumerate}
\item
    For every thin class $X$ there exists a strong well-ordering of $X$ (a well-ordering is strong
    iff every sub\emph{class} of X has a least element).
\item
    Every infinite class contains a countable infinite subclass.
\end{enumerate}
\end{proposition}
\begin{proof}
If an infinite class $X$ is thin then one can build a strong well-ordering
of $X$, applying Proposition~\ref{powerclass-codable} and \ax{Choice},
very much like the way one gets a well-ordering of a set using the axiom of
choice.
The least infinite initial segment of $X$
under that ordering will be countable.

If $X$ is not thin it contains an infinite subset $x$. 
Hence, there exists a bijective mapping $h$ from an infinite natural number onto
$x$. The class $h"\snat$ will be countable and infinite.
\end{proof}
\begin{remark}
The statement converse to Proposition~\ref{thin-wo},1) is true if we accept 
an additional axiom analogous to the axiom of chromatic classes
of \NCT.
\end{remark}

The following proposition describes small sets
in a way similar to the classical Dedekind's characterization of finite
sets (a set is finite iff it is not of equal cardinality with any
its proper subset).

\begin{proposition}[\THS] \label{small-Dedekind}
$\All X \Fbigl( \smal X\liff$
\[
    \All Y\ss X \Fbigl[Y\neq
    X \impl \All F:X\to Y \Fl(\text{"$F$ is not injective"}\Fr)\Fbigr]\Fbigr).
\]
\end{proposition}


\section{$\in$-structures and Zermelo universes}

\newcommand{\fcode}[1]{\lceil#1\rceil}
\newcommand{\True}{\opname{True}}

It is known (see~\cite{HKK,HK}) that saturation principles allow to
simulate structures satisfying the axioms of Zermelo set theory or even
{\ZFC} within nonstandard models of arithmetic.
From the other hand, in "fully saturated" nonstandard set theories (such as
E.~Nelson's Internal Set theory, {\NCT} or Hrbacek Set Theory of V.~Kanovei
and M.~Reeken~\cite{KR}) every $\in$--structure of standard size is isomorphic to an
$\in$--substructure of some hereditarily hyperfinite set.

In accordance with the above mentioned facts, the main result of this
section states that every thin semiset can be embedded in a thin subuniverse 
that satisfies axioms of Zermelo set theory with choice, subclasses in the sense of
{\THS} corresponding to subsets in the sense of the Zermelo subuniverse.

\begin{theorem}[\THS]\label{zermelo-universes}
For any thin class $X$ there exists a thin class $Z$ such that
\begin{enumerate}
\item
    $X=Z\cap x$ for some $x\in Z$;
\item
    $\All x\in Z \all C\ss x\cap Z \ex q\in Z \Fl( q \cap Z = C
    \Fr)$;
\item
    $\All x \Fl( x\ss Z \impl x\in Z \Fr)$;
\item \label{zc-axioms-true}
    All axioms of \ZCm{} (the Zermelo theory with the axiom
    of choice and without the axiom of regularity) are true in $Z$.
\end{enumerate}
\end{theorem}

\FACT
A class satisfying the conditions (2), (3) and (4) of
Theorem~\ref{zermelo-universes} is called a \em{Zermelo universe}.
\EFACT

\begin{remark}\label{ast-zermelo}
The definition of Zermelo universe and Theorem~\ref{zermelo-universes} are
very close in formulation to {\ZF}--classes and Cantorian axioms in \AST
(given in chapter 12 of \cite{VopS}). But the important condition (2) does not hold in AST.
\end{remark}

This theorem becomes a theorem of {\THS} if we give a formal meaning to
\itemref{zc-axioms-true} using encoding of formulas.

We fix some explicit coding, by standard sets, for symbols of logical
connectives, quantifiers, membership relation, 
punctuation signs and a countable set of variables, a coding 
for sets as parameters. \defterm{Formal formulas} are naturally defined within {THS}
by induction
as special (well-formed) sequences of codes. Every formula $\phi$ of \THS{}
gets its formal counterpart $\fcode{\phi}$ --- the code of $\phi$. Any
formal formula of small length with standard parameters is standard (as a
set).

The language $\SL(P)$ is defined as the class of all formal
formulas of small length with parameters from the class $P$. 

Evidently, the language $\SL(P)$ is thin for any thin $P$.

\begin{proposition}[\TFS]\label{true-formulas}
For any class $X$ there exists a unique class $T$ which
consists of closed formulas of $\SL(X)$ and satisfies the following
properties:
\begin{enumerate}
  \item $\fcode{x_1=x_2}\in T \liff x_1,x_2\in X \land x_1=x_2$;
  \item $\fcode{x_1 \in x_2}\in T \liff x_1,x_2\in X \land x_1\in x_2$;
  \item $\theta_1 \fcode{\lor} \theta_2 \in T \liff 
            \theta_1 \in T \lor \theta_2 \in T$;
  \item $\fcode{\lnot}\theta_1 \in T \liff \theta_1 \notin T$;
  \item $\fcode{\Ex} v \theta \in T \liff 
            \Ex x\in X \Fl( \theta_{v\to x} \in T \Fr)$,
\end{enumerate}
where $\theta_1,\theta_2$ are closed formulas of $\SL(X)$, $\theta$ is a
formula of  $\SL(X)$ with the only free (symbol of) variable $v$ and $
\theta_{v\to x}$ is obtained from $\theta$ by replacing $v$ with the code of
the set $x$.
\end{proposition}

We denote as $\True(X)$ the class, the existence of which is stated by
Proposition~\ref{true-formulas}.

For any closed formula $\phi$ with parameters from some class $X$
it is provable in \THS{} that
\[
    \phi^X \liff \fcode{\phi} \in \True(X),
\]
where $\phi^X$ is a relativization of $\phi$ to the class $X$.

We denote for any $\theta \in \SL(X)$
\[
    X \models_f\theta \iffdef \theta \in \True(X).
\]

Now we formalize (4) from Theorem~\ref{zermelo-universes}:
\[
\text{(4)\qquad} Z\models_f \theta \text{ for each $\theta$ such that "$\theta$ is an
axiom of {\ZCm}",}
\]
where the phrase in quotes is appropriately expressed as a formula of \THS.

\newcommand*{\Sd}{\opname{Sd}}
\newcommand*{\Def}{\opname{Def}}

We define the class $\Def(X)$ of sets
definable with a formula from $\SL(X)$ as follows:
\begin{gather*}
\Def(X) \eqdef \big\{ x: x = \coll{y: \bH\models \theta(y)}: \theta\in\SL(X) \text{
has exactly one free variable} \big\}\\
\end{gather*}

Obviously, $\Def(X)$ is a class. If $X$ is thin, $\Def(X)$ is also thin.

We will say that a class $C$ is an \emph{f--elementary submodel} of a class
$M\supseteq C$ 
(notation: $C\preccurlyeq_f M$) iff
\[
    C \models_f \phi \liff M \models_f \phi
\]
for any $\phi\in\SL(C)$.

\begin{theorem}[\TFS]\label{def-elem}
For any class $X$ the class $\Def(X)$ is an f--elementary submodel of $\bH$.
\end{theorem}
\begin{proof}
Note that $\K^{-1}(\min\coll{\K(a):\theta(a)})\in\Def(X)$ for $\theta\in\SL(X)$.
\end{proof}
\begin{corollary}
In the theory \THS{} without the transfer axioms, the transfer axioms follow
from the statement $\Def(\emp)=\stan$.
\end{corollary}

\begin{proof}[Proof of Theorem \ref{zermelo-universes}]
Using the axiom of dependent choices we construct the sequence of structures
$S_n$ as follows.

We start from some thin class $S_0 \supseteq X$ such that $S_0 \preccurlyeq \bH$.
Such a class does exist by Theorem~\ref{def-elem}.

Given a class $S_n$, using the axioms of exponentiation and choice, 
we can choose a thin class $S_{n+1}$ to satisfy the following properties:
\begin{enumerate}
\item
    $\All C\ss S_n \ex y\in S_{n+1} \Fbigl( y\cap S_n = C \land \Fl(
    \Set(C)\impl y=C \Fr) \Fbigr)$;
\item
    $\All x,y\in S_{n+1} \Fl( x \neq y \impl x\cap S_n \neq y \cap S_n \Fr)$;
\item
    $\All x\in S_{n+1}\setminus S_n \all y\in S_n \Fl( x\notin y \Fr )$.
\end{enumerate}

We put
\[
    Z = \bigcup_{n\in\snat} S_n.
\]

It is easy to see that the axioms of extensionality, union, power set,
separation, infinity and choice hold in $Z$.

\end{proof}

\begin{corollary}
\begin{enumerate}
    \item
        \THS{} is not a conservative extension of \ZFfin;
    \item
        {\THS} is strictly stronger than \ZCm.
\end{enumerate}
\end{corollary}

If we take $X=\stan$ in the conditions of the previous Theorem and apply
transfer, we get immediately the following theorem.

\begin{theorem} \label{finite-rich}
Every statement about finite sets provable in \ZCm{} holds in \THS{} as
well.
\end{theorem}

\begin{remark}
Non-standard extensions of superstructures over a thin class can be
constructed easily in {\THS} as thin $\in$--structures which allows to use
"essentially external" methods of nonstandard analysis such as nonstandard
hulls of Banach spaces or Loeb measures.
\end{remark}


\section{Real numbers}

Real numbers can be introduced in {\THS} in a quite usual and straightforward
way --- as elements of a complete linearly ordered field.

We introduce explicitly the rational numbers first.

\begin{fakt}
First of all define operations $+$ and $\cdot$ on $\N$ by the
formulas
$$
x+y=\sharp (x\cup\{0\}\times y);\ x\cdot y=\sharp(x\times y)
$$

Obviously the introduced operations satisfy the the classical
recursive definitions of addition and multiplication of natural
numbers. It is easy to see that the subclass $\snat$ of $\N$ is
closed under these operations.
\end{fakt}

\begin{fakt}
Usually the ring of integers is defined as quotient set of
$\N\times\N$ under the equivalence relation
$$
\<a,b\>\sim\<a_1,b_1\>\Liff a+b_1=a_1+b.
$$

Since in our case $\N$ is a class we must define the quotient
class by a system of representatives.

Thus, the class $\Z$ of integers can be defined e.g. by the
formula
$$
\Z=\{0\}\times\N\cup\N\times\{0\}
$$
with obviously defined addition, multiplication and linear order
relation. It is easy to prove also that the thin class $\sint$ of
standard elements of $\Z$ is a subring of $\Z$.
\end{fakt}

\begin{fakt}
The field $\Q$ is defined as the quotient field of the integral
domain $\Z$. As before we must define this quotient field by a
system of representatives, e.g. by the formula
$$
\Q=\{\< a,b\>\in\Z\times\Z\ |\ b\neq 0,\ \opname{gcd}(a,b)=1\}
$$

Once again it is easy to prove that the thin class $\srat$ of
standard elements of $\Q$ is a subfield of $\Q$.
\end{fakt}

\begin{fakt} A class $\<R;{+},\cdot,\le\>$ is called a \defterm{field
of real numbers}\ iff it satisfies the axioms of linearly ordered
field and the following completeness property:
\[
    \text{every bounded above \emph{subclass} of $R$ has a supremum.}
\]
\end{fakt}
\begin{theorem} \label{Reals}(\THS)
\begin{enumerate}
\item There exists a thin class that is a field of real numbers.

\item The field $\srat$ is dense in a field of real numbers.

\item  Any two fields of real numbers are isomorphic.

\end{enumerate}
\end{theorem}
\begin{proof}[Proof sketch]
The proof quite repeats the classical one. We can choose any usual way of constructing real numbers. Take, for example,
Dedekind cuts. Due to axiom of exponentiation $\coll{C: C\ss \srat} = \coll{
P"\coll{c} : c\in\dom(P) }$ for some thin class $P$. Every Dedekind cut can
be identified then with an element of $\dom(P)$ and we build the field of
real numbers as a subclass of $\dom(P)$.

The classical proofs of (2) and (3) can also be transferred easily to {\THS}
(see also Theorem~\ref{thsint}).
\end{proof}

\begin{remark}
In every Zermelo subuniverse $Z$, the field of reals in the sense of $Z$ is
a field of reals in the global sense.
\end{remark}

\begin{fakt}In what follows we fix some field of real numbers
$\<\sreal; {+},{\cdot},{\leq}\>$, and call it \emph{the} field of reals.
\end{fakt}

\begin{remark}
There is no definable field of real numbers in {\THS} (see
Proposition~\ref{no-definable-uncountable}).
\end{remark}

As in non-standard analysis, every bounded rational number has a standard
part in $\sreal$.

Put $\Q_b=\coll{ x\in\Q : \Ex r\in\srat\,(|x|<r)}$. We call elements
of $\Q_b$ bounded rationals. Obviously $\Q_b$ is a subring of
$\Q$ and $\srat\ss\Q_b$.

Let $\mu(0)=\coll{\a\in\Q_b : \All r\in\srat \Fl(r>0\impl |\a|<r\Fr)}$.
Then $\mu(0)\ss\Q_b$ is an ideal in $\Q_b$.

\begin{theorem}\label{stpart}
There exists a unique surjective homomorphism $\st:\Q_b\to\R$.
The kernel $\ker(\st)=\mu(0)$.
\end{theorem}

The real number $\st(x)$ is called \emph{the standard part of a
bounded rational $x$}.


\section{Ordinary mathematics in \THS}

Intuitively, thin classes behave exactly as usual infinite sets. We would
like to transfer notions and results of ordinary mathematics to systems of
thin classes. The informal principle is:
\begin{itemize}
\item[] \emph{Everything that is true in ordinary mathematics about sets, 
their subsets, powersets and so on is true in {\THS} about thin classes, 
their subclasses, systems of their subclasses and so on.}
\end{itemize}

Note that cartesian products are implemented within iterated powersets; 
finiteness can be expressed as Dededkind finiteness and is equivalent to
smallness by Proposition~\ref{small-Dedekind}.

In what follows we will give a formal account of the formulated principle.

\begin{fakt}\label{topology-example} 
As an example, we would like to say whether a system
$\mathcal{T}$ of subclasses of a thin class $X$ is a topology on $X$. If
$\mathcal{T}=\coll{T"\coll{d}: d\in D}$ this can be expressed in the
following way:
\begin{multline}\label{m-topol-space}
    T"D = X \land
    \All d_1,d_2\in D \Fl( T"\coll{d_1}\cap T"\coll{d_2}\neq \emp \Fr)
\land\\
    \All D'\ss D \Ex d \Fl( \bigcup_{e\in D'} T"\coll{e} = T"\coll{d} \Fr).
\end{multline}

If $X$ is represented in a Zermelo universe $Z$ by an element $x\in Z$ ( $Z\cap x =
X$) then $\mathcal{T}$ is also
represented in $Z$ by some $t\in Z$: 
\[
    \all Y \Fbigl( 
        Y\in\mathcal{T} \liff \Ex y\in t\cap Z \Fl( y\cap Z = Y \Fr)
    \Fbigr),
\]
and (\ref{m-topol-space}) is true iff $Z\models \text{"$t$ is a topology on
$x$"}$.
\end{fakt}

\FACT
In a more generic setting we may need to refer to higher levels of cumulative
hierarchy over some thin class. Some encoding is necessary for that. 
To describe a general situation and abstract from a particular encoding
of systems of classes (as we did in~\ref{topology-example})
we consider extensional systems over thin classes.
\EFACT

\newcommand{\cE}{\operatorname{\mathcal{E}}}
\begin{fakt}
A system of classes (see~\ref{system-of-classes}) $\cX$ equipped by a system of pairs of classes ${\cE}$ is
called a (thin) \defterm{extensional system} over a thin class $A$ iff the
following conditions hold:
\begin{enumerate}
\item
    $X \cE Y \impl X\in \cX \land Y\in \cX$;
\item
    $A\in \cX$, $A_{\cE} = A$, $\All a\in A\Fl( a_{\cE} = a\cap A \Fr)$, 
    where $Y_{\cE} \eqdef \coll{Z : Z\cE Y}$;
\item
    $\All X,Y\in \cX\setminus A \Fbigl( \All Z \Fl( Z \cE X \liff Z \cE Y
        \Fr) \impl X=Y \Fbigr)$;
\item
    $\All X\in\cX \ex\qsuper{\smal} k \Fl(
        \underbrace{\cup_{\cE} \dots \cup_{\cE}}_{\text{$k$ times}} X_{\cE}\ss A$
    \Fr) \\
    where $\cup_{\cE} \mathcal{S} \eqdef \coll{ Y: \ex Z
                    \Fl( Z\in \mathcal{S} \land Y\cE Z\Fr)}$.
\end{enumerate}

We put \\
$\cX_0 = A$; $\cX_1=\coll{X:X_{\cE}\ss A}$; 
$\cX_k = \coll{X\in\cX:
    \underbrace{\cup_{\cE} \dots \cup_{\cE}}_{\text{$k-1$ times}} 
    X_{\cE}\ss A}$, $k>1$.

A system $\cX$ is called \defterm{$k$--full} iff
\[ \All \cY\ss\cX_k \ex Y\in\cX \Fl( Y_{\cE}=\cY \Fr).\]
\end{fakt}

\begin{fakt}
We define a \defterm{formula of ordinary mathematics} (o.m.--formula) to be an
$\in$--formula $\phi(A,X_1,\ldots,X_n)$ where all quantifiers 
have the form $\Ex x\in \cP^k(A)$ or $\All x\in \cP^k(A)$ where 
each quantifier has its own natural
number $k$ of iterations of the powerset operation $\cP$. 
(Formally, "$\All x\in \cP(A) \dots$" is to be read as 
$\All x \Fbigl( \All z \Fl( z\in x \impl z\in A \Fr) \impl \dots \Fbigr)$,
and so on).
The maximal number of iterations of $\cP$ in the bounding terms 
of $\phi$ is called the \defterm{height} of $\phi$.

For any o.m.--formula $\phi$ of height $k$ the truth of
$\phi(A,X_1,\ldots,X_n)$ in a $k$-full extensional system
$\cX$ over $A$ ($X_i\in\cX$) is defined in a straightforward way (we omit the obvious details). We write
$\cX\models\phi$ if $\phi$ is true in $\cX$.
\end{fakt}

\newcommand{\cJ}{\mathcal{J}}
\begin{theorem}\label{om-true}
Let $\cX$ be an extensional system over a thin class $A$. 
    Suppose $A$ is represented in a Zermelo universe  $Z$: $A=Z\cap a$ for some $a\in Z$.
    Then there exists a \emph{unique} embedding 
    $\cJ : \cX \to Z$ such that
    \[
        \cJ(A)=a \land \All X, Y \in \cX \Fl( X \cE Y \liff \cJ(X)\in\cJ(Y) \Fr).
    \]
    Moreover, if $\cX$ is $k$--full then for any $X_1,\ldots,X_n \in \cX$ and any o.m.--formula
    $\phi(A,X_1,\ldots,X_n)$ of height $\le k$ we have
    \[
        \cX\models\phi(A,X_1,\ldots,X_n) \liff
        Z\models\phi(\cJ(A),\cJ(X_1),\ldots,\cJ(X_n)).
    \]
\end{theorem}

The proof is quite straightforward.

\begin{fakt}
Theorem \ref{om-true} allows to extend the definition of truth of
o.m.-formulas so that it can be applied to extensional systems
that not necessarily are full to the height of the formula.

Let $\cX$ be an
extensional system over a thin class $A$, $X_1,\ldots,X_n\in\cX$ and
$\phi(A,X_1,\ldots,X_n)$ be an o.m.-formula of height $k$. Then we say that
\linebreak $\phi(A,X_1,\ldots,X_n)$ is true for $\cX$ iff
$\cX'\models\phi(A,X_1,\ldots,X_n)$ for some (and then for any) $k$--full
extensional system $\cX'\supseteq \cX$ over $A$.
\end{fakt}

\begin{theorem}\label{thsint}
Suppose $\All A\all X_1\in\cP^{k_1}\dots \all X_n\in\cP^{k_n} \f
\phi(A,X_1,\ldots,X_n)$ is a theorem of ordinary mathematics provable in
{\ZCm}. Let $X_1\in\cX_{k_1},\ldots,X_n\in\cX_{k_n}$ in an extensional
system $\cX$ over a thin class $A$. Then $\phi(A,X_1,\ldots,X_n)$ is true
for $\cX$.
\end{theorem}

It is easy to see, for example, that the statement (3) of Theorem~\ref{Reals}
follows from this theorem.

\section{Interpretaions of \THS}

\begin{theorem}\ \label{interps}
    The collection of all subclasses of the
    set $V_\omega$ of hereditarily finite sets together with the original 
    membership relation gives an interpretation of \THS{} in \NCT. Moreover,
    under this interpretation
    \begin{enumerate}
    \item
        sets are exactly finite subsets of $V_\omega$;
    \item
        thin classes are exactly subsemisets of $V_\omega$ of standard size.
    \end{enumerate}
\end{theorem}
\begin{proof}
By Theorem~3.12 of \cite{AG} a set is S-finite in {\NCT} (that is having a
standard finite cardinality) iff every subclass of it is a set. Corollary~4.12
 of \cite{AG} states that a semiset $X$ has a standard size iff every
subset of $X$ is S-finite. So we have that small sets are interpreted as
S-finite, and thin classes are interpreted as semisets of standard size.
\ax{Exten} holds obviously. Axiom \ax{Class} follows from
Corollary~4.17 of \cite{AG} stating that any formula in which only semisets are
quantified is equivalent to a normal formula. Axioms \ax{Set} and \ax{Ind}
can be derived from the theorem of Sochor mentioned in the proof of
Theorem~\ref{zffin}.
Proposition~4.13 of \cite{AG} says that any semiset of stanard size can be
embedded into a set of any given infinitely large cardinailty. Axiom of thin
semisets follows. The truth of \ax{Comp} can be derived easily from the
Saturation Theorem~4.7 of \cite{AG}.

The Choice Theorem~4.20 of \cite{AG} implies \ax{Choice}. The axiom of
exponentiation \ax{Exp} also follows from the Choice theorem because if
$\kappa$
is the "standard size" of a semiset $X$ then $2^\kappa$ is the "standard size"
of the class of its subclasses.

Theorem 1.10 of Kanovei and Reeken~\cite{KR1} shows that the scheme of
dependent choices for sets holds in {\BST} (Bounded Set Theory). Since
{\NCT} is a conservative extension of {\BST} (\cite[Theorem 5.1]{AG}) and 
semisets are uniformly parameterized by sets (\cite[Theorem 4.16]{AG}),
the scheme \ax{DC} is also true.
\end{proof}

\begin{remark}
In a model of E.~Nelson's~\IST{}\cite{Ne1} the collection of all subclasses of $V_\omega$ in
that model gives a model of \THSZ{} but does not give a model of the full \THS. The reason is that there exists a
subclass $O$ of $V_\omega$ that can be one-to-one mapped onto the class of all
standard sets and therefore will be thin but neither \ax{Choice} nor \ax{Exp} can be proved for $O$.
\end{remark}

\begin{proposition}\label{no-definable-uncountable}
There is no formula $\Phi$ with one free variable such that $\Ex! X \f
\Phi(X) \land \All X \Fl( \Phi(X) \impl \text{"$X$ is an uncountable thin
class" \Fr)}$ would be a theorem of {\THS}.
\end{proposition}
\begin{proof}
Under the interpretation of {\THS} in {\NCT} described above every formula
of {\THS} gets translated to a normal formula of {\NCT}. By Proposition~7
from~\cite{AH} any class of standard size defined by a formula without
parameters consists of standard elements\footnote{The proof is given in
\cite{AH} for
{\BST} but can be transfered literally to {\NCT}.}. Therefore any thin class
definable by a formula in such a model has to consist of hereditarily small
sets, and hence cannot be uncountable.
\end{proof}

\newcommand*{\ADef}{\opname{ADef}}

For any class $X$, denote
\[
    \ADef(X) = 
    \bigcap_{n\in \Def(X) \setminus \stan} \coll{ x: \K(x) < n }.
\]

Due to compactness, $\ADef(X)$ is nonempty for any thin class $X\neq\stan$.

\begin{proposition}
For any thin class $X\neq\stan$, the collection of all subclasses of the class $\ADef(X)$ 
together with the original membership relation forms an interpretation of
\THS{} in \THS.
\end{proposition}

Thus, there are interpretations of \THS{} in \THS{} of any size: $\ADef({l})$ is a
subclass of a set having less than $l$ elements. Moreover, the universe
$\bH$ of all sets can be thought of as a subclass of some highly unfeasibly
large hyperfinite set.

\section{Indiscernibilty equivalences and locally compact topological spaces}
\label{indiscernibility}

In this section we discuss how the approach to continuous
structures that considers them as the images of accessible parts
of certain hyperfinite structures under identifying indiscernible elements
can be formalized in \THS.

\newcommand{\eps}{\operatorname{\stackrel{E}{\approx}}}
Let $x$ be a set and $\eps$ be a $\pi$--equivalence relation on
$x$ which means that $\eps = \cap E$ where $E$ is a thin class of
subsets of $x\times x$.

\begin{fakt}
Let $X\ss x$ be a $\s$-subset of $X$, which means that $X$ is the
union of a thin class of subsets of $X$. The relation $\eps$ is
called an \defterm{indiscernibility equivalence} on $X$ iff
$$
    \Allinf u\ss X \ex a,b \in u \ ( a\neq b \land a \eps b ).
$$
\end{fakt}

\begin{fakt}[Example]
To consider an example of an indiscernibility relation we introduce
the following notation.

Let $\a,\b\in\Q$. Then
\begin{enumerate}
\item $\a\approx\b\iffdef \a-\b\in\mu(0)$ (cf. Theorem
\ref{stpart});

\item $\a\sim\infty\iffdef\a\in\Q\setminus\Q_b$.
\end{enumerate}

Note that $\N\ni n\sim\infty\Liff n\in\N\setminus\snat$.

Fix $n,m\in\N\setminus\snat$ such that $\frac
nm\sim\infty$. Let $x=\{\pm\frac kl\ | 0\leq k\leq n,\ 0<l\leq
m\}$. Consider the restriction of the relation $\approx$ to $x$.
We denote this restriction also by $\approx$ in this example. Let
$X=x\cap\Q_b$. Obviously $X$ is a $\s$-class. Let us show that
$\approx$ is an indiscernibility relation on $X$~\footnote{This is
obvious for those, who are familiar with nonstandard analysis}.
Indeed, let $u\ss X$. Put $v=\{|a-b|\ | \ a,b\in u,\ a\neq b\}$.
If $\min v\approx 0$ then for every $n\in\snat$ the set$u_n=\{
a,b\in u\ |\  a\neq b,\ |a-b|<\frac 1n\}\neq\emptyset$. Since the
decreasing countable sequence $\{u_n\ |\ n\in\snat\}$ consists of
nonempty sets, its intersection is also nonempty by 
Proposition~\ref{saturation-extension}. Thus there exist $a,b\in u$ such that
$a\neq b$ but $a\approx b$. Therefore, in this case our statement
is proved.

Let now $\min v>\d>0$, where $\d\in\srat$. This implies that the
map $\st:\Q_b\impl\R$ defined in Theorem \ref{stpart} is injective
on $u$ and $\min\{|\st(a)-\st(b)|\ |\ a\neq b,\ a,b\in U\}\geq\d$

On the other hand the set $\st"u$ is bounded. Indeed, since
$u\ss\Q_b$, we have \newline $st"u\ss[\st(\min u)-1,\st\max u+1]$.
Applying the theorem of ordinary mathematics, which states that
every infinite bounded set of reals has an accumulation point we
obtain, that $st'u$ is finite and, thus, $u$ is small.

Using Theorem \ref{thsint} we can easily formalize this
consideration in \THS.
\end{fakt}

\begin{theorem} \label{ind_prop} The following conditions are equivalent:
\begin{enumerate}
\item
    $\eps$ is an indiscernibility equivalence;
\item
    $\All e\in E \All y\ss X \Ex\qsuper{\smal} v\ss X ( y \ss e"v )$;
\item
    $\Exthin N\ss X (\All a\in X \ex n\in N ( a \eps n )\land \All n,m\in N (n\eps m\impl n=m))$.
\end{enumerate}
\end{theorem}
\begin{proof}
$(1)\impl (2)$. If (2) does not hold then there exist $e\in E$ and $y\ss X$ such
that $a\setminus(y"v)\neq\emp$ for any small $v\ss y$. 
Using induction for $\snat$, we will
find for any small $n$ an injective function from $n$ into $y$ such that its
range consists of pairwise $e$-non-equivalent elements. 
Applying prolongation, we get an infinite set of pairwise $e$-non-equivalent elements.

$(2)\impl(3)$ Assume $X=\cup D$. Using choice and prolongation we can assign to any $e\in E$
and $d\in D$ a set $v_{ed}$ such that $d\ss e"v_{ed}$. Take any $d\in D$ and
elements $x,y\in d$. Suppose $\All e\in E\ex z\in v_{ed} \Fl( x\in
e"\coll{z} \land y\in e"\coll{z} \Fr)$. Then $x\eps y$. Indeed, take any
$e\in E$. Since $E$ is an equivalence relation there is an $e_1\in E$ such
that $e_1\circ e_1^{-1}\ss e$. Since, for some $z$, $\<x,z\>\in e_1 \land
\<y,z\>\in e_1$, we have $\<x,y\>\in e$.

Now, since the class $B=\coll{e"\coll{z} : z\in v_{ed} \land e\in E\land d\in D}$ is thin, 
we can choose $N$ in such a way that $N$ contains exactly one element in the
intersection of each centered subfamily of $B$ ( $C\ss B$ is centered iff
any small subset of $B$ has a non-empty intersection). It can be checked
easily that $N$ is as required in (3).

The implication $(3)\impl (1)$ is straightforward.
\end{proof}

Obviously, a class $N$ satisfying the condition 3 of Theorem
\ref{ind_prop}, generally speaking, is not unique. Fix any such
$N$.  \newcommand{\inter}[1]{\stackrel{\circ}{#1}} For $u\ss X$
define $\inter u=\{a\in u : \All b\in X(a\eps b\impl b\in u)\}$
and $\inter u^{\#}=\inter u \cap N$.

\begin{proposition}[\THS] \label{top}
\begin{multline}
\Fbigl[\All n\in N\ex u\ss X(n\in\inter u^{\#})\Fbigr] \land 
\All n\in N\ \All u_1,u_2\ss X \\
    \Fbigl( n\in\inter u_1^{\#}\cap\inter u_2^{\#}\impl
    \ex u\ss X \Fl(n\in \inter u^{\#}\ss\inter u_1^{\#}\cap\inter u_2^{\#}\Fr)
    \Fbigr)
\end{multline}
\end{proposition}

\newcommand{\T}{\mathcal{T}}
Semantically, this proposition means that the system of classes
$\T=\{\inter u^{\#}\ |\ u\ss X\}$ forms a base of topology on $N$.

It can be proved that this topology is locally compact. Let us
show how this statement can be formulated explicitly in \THS. Let $C\ss N$.
We say that $F$ is an open covering of $C$ if $F$ is a function,
$\dom F=I$ is a thin class and $\All i\in I F(i)\ss X$ and
$N\ss\bigcup_{i\in I} \inter F(i)^{\#}$. A class $C\ss N$ is
compact iff $\All F(F\ \mbox{is an open covering of}\
C)\longrightarrow \Ex\qsuper{\smal} p\ss\dom F(C\ss\bigcup_{i\in p}\inter
F(i)^{\#})$ (cf. Remark \ref{small-Dedekind}). Now, "$(N,\T)$ is a
locally compact space" is equivalent to the following \THS-formula:
$$
\All n\in N\Ex u\ss X(n\in u^{\#}\land\Ex C\ss N((C\ \mbox{is
compact})\land (u^{\#}\ss C)). \eqno(LC)
$$

Similar approach to the construction of locally compact spaces was
developed in \cite{Gor} in terms of nonstandard analysis for the
case of locally compact abelian groups. The proof of Theorem 2.2.4
of \cite{Gor} can be easily transformed to a proof of
(LC) in \THS.

It is easy to see that any locally compact space can be
represented as a quotient class $N$ constructed by an appropriate
triple $\<x, X,\eps\>$, where $X$ is a thin subclass of a set $x$
and $\eps$ is an indiscernibility relation on $X$.

Let us consider such representations of the field $\R$ in more detail. As it was
mentioned in the Introduction they can be considered as numerical
systems that simulate reals in an idealized computer of infinite
(i.e. hyperfinite non-small) memory. This point of view gives a
motivation for the following definition.

\begin{fakt} \label{compar}
Consider a tuple $R=(\< r;\oplus,\odot\>; R_b,\rho)$, where 
$r$ is a hyperfinite set with binary operations
$\oplus,\odot$ on it. We say that $R$ is a
\defterm{hyperfinite computer arithmetic} if $R_b\ss r$ is a $\s$-class and
$\rho$ is a $\pi$-equivalence relation  on $r$ such that
\begin{enumerate}
\item $\rho$ is an indiscernibility relation on $R_b$;

\item $R_b$ is closed under the operations $\oplus$ and $\odot$;

\item $\rho$ is a congruence relation on $R_b$;

\item the quotient algebra $R_b/\rho$ is topologically isomorphic
to the field $\R$.
\end{enumerate}
\end{fakt}

The previous considerations show that this definition can be
formalized in \THS.

We interpret the elements of $R_b$ the same way as it was
discussed in the Introduction. Elements of $R_b$ are computer
reals that are not too big, i.e. not too close to the boundary of
the computer memory. Obviously, this is not a definition in the
framework of the classical mathematics. That is why $R_b$ is a
proper semiset. Operating with these numbers does not imply the
overfilling of the memory. Thus the computer operations restricted
to $R_b$ approximate the corresponding operations on reals. This
fact is formalized in the statement 3 of definition~\ref{compar}.

\newcommand{\om}{\omega}

{\bf Example 2}. Let $0<\e\in\mu(0)$ and $\om\in\N\setminus\snat$
be such that $\om\e\sim\infty$. Consider a tuple 
$R(\om,\e)=(\<r_{\om};\oplus,\odot\>; R_b,\rho)$, where
\begin{enumerate}
\item $r_{\om}=\{-\om,\dots,\om\}$;

\item the operation $\oplus$ is the addition modulo $2\om+1$;

\item the operation $\odot$ is defined by the formula
$$
k\odot m =[km\e](\mod (2\om+1)),
$$
where $k,m\in r_{\om}$ and $[\a]$ is the integral part of a real
number $\a$;

\item $R_b=\{k\in r_{\om}\ |\ k\e\in\Q_b\}$;

\item an equivalence relation $\rho\ss r_{\om}\times r_{\om}$ is
such that
$$
k\rho m\Liff k\e\approx m\e
$$
\end{enumerate}

\begin{proposition} \label{example}
The tuple $R(\om,\e)$ is a nonstandard computer arithmetic
\end{proposition}

$\rhd$ It is easy to see that if $k,m\in\R_b$ then $k\oplus m=k+m$
and $k\odot m=n$, where $n\e\leq km\e^2<(n+1)\e$.

Define the map $F: R_b\impl\R$ by the formula $F(k)=\st(k\e)$. It
is easy to see that for all $k,m\in R_b$ holds
\begin{enumerate}
\item $F(k)=F(m)\Liff k\rho m$;

\item $\All q\in\Q_b\Ex k\in R_b (q\approx k\e)$;

\item $(k\odot m)\e\approx k\e\cdot m\e$.
\end{enumerate}

These properties prove that $F$ is a surjective homomorphism and
that $R_b/\rho$ is isomorphic to $\R$ $\lhd$

The nonstandard computer arithmetic $R(\om,\e)$ discussed in
Example 2 is not a hyperfinite version of the computer arithmetic,
which is implemented in existing computers. The last one is based
on the floating point representation of reals. We will call it
FP-arithmetic. Its hyperfinite version was discussed in \cite{ggh}
in terms of nonstandard analysis. The computer arithmetic
$R(\om,\e)$ (no matter standard or nonstandard $\om$ and $\e$ are
considered) has some better than FP-arithmetic algebraic
properties. Indeed, it is well-known that the addition and the
multiplication in FP-arithmetic are neither associative, nor
distributive, while $\<r_,\oplus\>$ is an abelian group. However,
multiplication in $R(\om,\e)$ is not associative and the law of
distributivity also fails \cite{ggh}.

It is not quite clear how the good algebraic properties of
numerical systems would affect on the quality of numerical
computations. It was shown \cite{Gor} that the convergence
properties of approximation of the Fourier Transformation on $\R$
by sampling of its kernel are better when the result of this
sampling is the matrix of Finite Fourier Transformation, i.e. when
we approximate the additive group $\R$ is by finite abelian
groups. The theory of approximation of locally compact groups by
finite abelian groups was developed in \cite{Gor}. It can be
proved (cf. \cite{ggh}, where similar questions where discussed in
terms of nonstandard analysis) that the problem of approximation
of locally compact algebraic systems by finite ones can be reduce
to a problem of representation of locally compact systems by the
quotients of $\s$-subsystems of hyperfinite systems under
indiscernibility relations.

The following theorem demonstrates the restrictions that occur on
the way of construction of computer arithmetics with the best
possible algebraic properties.

\begin{theorem} \label{nonring}
There does not exist a hyperfinite computer arithmetic 
$R=(\<r;\oplus,\odot\>; R_b,\rho)$ such that $\<r;\oplus,\odot\>$ is
an associative ring (even non-commutative).
\end{theorem}

A similar theorem about finite approximations of locally compact
fields was proved in \cite{ggh} (see also \cite{gg}). Theorem
\ref{nonring} is a little bit more general than those of
\cite{ggh}. However, the proof presented in \cite{ggh} can be
easily adjusted to Theorem \ref{nonring}. This proof can be
formalized in \THS.

\newcommand{\Allb}{\All^b}
\newcommand{\Exb}{\Ex^b}

Let $L$ be the first order language in the signature 
$\s_1=\<+,\cdot\>$ and $L_h$ - the first order language in the signature
$\s_2=\<\oplus,\odot;R_b,\rho\>$. Here $\oplus$ and $\odot$ are
symbols of binary operations, $R_b$ is a symbol of a unary
predicate and $\rho$ - a symbol of binary predicate.

Formulas of $L_h$ have the natural interpretation in any
hyperfinite computer arithmetic $R=(\<r;\oplus,\odot\>;
R_b,\rho)$. We use notations $\Allb$ and $\Exb$ for the universal
and the existential quantifiers restricted to $R_b$.

Let $t$ be a term in the signature $\s_1$. Replace each occurrence
of the operation $+$ ($\cdot$) in $t$ by $\oplus$ ($\odot$). The
obtained term in the signature $\s_2$ will be denoted by $t_h$.

\renewcommand{\f}{\varphi}

Let $\f(x_1,\dots,x_n)$ be an $L$-formula. Denote by $\f_h$ the
$L_h$-formula obtained from $\f$ by replacement of each atomic
formula $t=s$ by $t_h\,\rho\,s_h$ and each quantifier $Qx$ by
$Q^bx$. Here $t,s$ are $\s_1$-terms. We call $\phi_h$ \emph{the
hyperfinite analog} of $\phi$.

Let $R=(\<r;\oplus,\odot\>; R_b,\rho)$ be a hyperfinite
computer arithmetic. By Definition \ref{compar} (4) there exists a
a homomorphism $\psi:R_b\impl\R$ such that $\psi(a)=\psi(b)\Liff
a\,\rho\,b$. This homomorphism may not be unique. We call $\psi$ a
canonical isomorphism. It may not be unique. The following theorem
follows immediately from the definitions.

\begin{theorem} \label{transfer}
For any $L$-formula $\f(x_1,\dots,x_n)$ the following statement
holds.

If $R=(\<r;\oplus,\odot\>; R_b,\rho)$ is a hyperfinite computer
arithmetic, $\psi:R_b\impl\R$, and $a_1,\dots a_n\in R_b$ then
$$
R\models\f_h(a_1,\dots, a_n)\Liff \R\models\f((\psi(a_1),\dots,\psi(a_n))
$$

\end{theorem}

This theorem gives some qualitative formalization of the fact that
if we operate with relatively small numbers, so that the memory
overfilling cannot occur during the computations, we obtain
results that are approximately true at least, when we deal with
algebraic statements that can be formalized in the language $L$.

A version of Theorem~\ref{transfer} can be formulated in the
language of classical mathematics only for some specific $L$-
formulas - the \emph{positive bounded formulas} \cite{ggh}.

The investigation of the correlation between continuous and
computer mathematics in terms of \THS{} (i.e. on the qualitative
level) for higher order properties (e.g. formulated in the
language of type theory) is an interesting problem. It can help to
discover some new phenomena concerning numerical investigation of
some more complicated structures.

\bibliographystyle{amsplain}

\end{document}